\documentclass{ifacconf}

\usepackage{xcolor}
\usepackage{graphicx}
\usepackage{xspace}
\usepackage{algorithm}
\usepackage{algpseudocode}
\algrenewcommand\algorithmicrequire{\textbf{Input:}}

\usepackage{amsmath,amsfonts,amssymb}
\usepackage{bm,fixmath} 
\usepackage{tikz}
\usetikzlibrary{arrows, positioning, calc}

\newtheorem{assumption}{Assumption}

\newtheorem{proposition}{Proposition}

\DeclareMathOperator*{\argmin}{arg\,min}

\newcommand{\N}{\mathbb{N}}
\newcommand{\R}{\mathbb{R}}

\newcommand{\norm}[1]{\left\lVert#1\right\rVert}

\newcommand{\x}{\mathbold{x}}

\newcommand{\w}{\mathbold{w}}

\newcommand{\A}{\mathbold{A}}
\newcommand{\B}{\mathbold{B}}

\newcommand{\G}{\mathbold{G}}

\usepackage{accents}
\newcommand{\lmax}{\bar{\lambda}}
\newcommand{\lmin}{\underline{\lambda}}

\usepackage{natbib} 
\begin{document}
\begin{frontmatter}

\title{Self-Identifying Internal Model-Based Online Optimization\thanksref{footnoteinfo}} 

\thanks[footnoteinfo]{
Corresponding author: N. Bastianello (\texttt{nicolba@kth.se})\\
The work of N. Bastianello was partially supported by the EU Horizon Research and Innovation Actions program under Grant 101070162.\\
This work has also been partially supported by the Wallenberg AI, Autonomous Systems and Software Program (WASP), funded by the Knut and Alice Wallenberg Foundation.}

\author[First]{Wouter J. A. van Weerelt} 
\author[First]{Lantian Zhang} 
\author[First]{Silun Zhang}
\author[Second]{Nicola Bastianello}

\address[First]{Department of Mathematics, KTH Royal Institute of Technology, Stockholm, Sweden}
\address[Second]{School of Electrical Engineering and Computer Science, and Digital Futures, KTH Royal Institute of Technology, Stockholm, Sweden}

\begin{abstract}
In this paper, we propose a novel online optimization algorithm built by combining ideas from 
control theory and system identification.
The foundation of our algorithm is a control-based design that makes use of the internal model of the online problem. 
Since such prior knowledge of this internal model might not be available in practice, we incorporate an identification routine that learns this model on the fly.
%
The algorithm is designed starting from quadratic online problems but can be applied to general problems.
For quadratic cases, we characterize the asymptotic convergence to the optimal solution trajectory.
We compare the proposed algorithm with existing approaches,
and demonstrate how the identification routine ensures its adaptability to changes in the underlying internal model. Numerical results also indicate strong performance beyond the quadratic setting.
\end{abstract}

\begin{keyword}
online optimization, online learning, system identification, online gradient descent, control-based optimization
\end{keyword}

\end{frontmatter}

\section{Introduction} 

The technological advances of recent years have increased the available sources of high resolution, streaming data in several applications, including the power grid, transportation networks, connected consumer devices \citep{simonetto,dall'anese}.
Leveraging these data at the relevant time-scales, \textit{e.g.} for control \citep{liao-mcpherson, paternain, chachuat}, signal processing, \citep{natali, hall, fosson}, machine learning \citep{shalev, dixit, rakhlin}, therefore requires the solution of \textit{online optimization} problems.

Online optimization problems are characterized by time-varying cost functions, which capture the evolving nature of the dynamic environments from which they arise \citep{simonetto}.
Formally, the problem of interest is
\begin{align}\label{eq:objective}
    \x^*_k =\argmin_{\x \in \mathbb{R}^n}f_k(\x),\hspace{2ex} k \in \mathbb{N},
\end{align}
where consecutive optimization instances arrive at intervals of length $T_s > 0$.


In addressing such problems, two categories of algorithms are typically considered: \textit{unstructured} and \textit{structured} \citep{simonetto}. By unstructured we refer to algorithms which solve the optimization problem revealed at each time step $k$, similar to traditional batch algorithms; online gradient descent is a typical example.
However, this approach does not tailor the design of the algorithms to the dynamic nature of the problem, thus in general we can only guarantee convergence to a neighborhood of the solution trajectory $\{ \x_k^* \}_{k \in \N}$ (assuming uniqueness to simplify this discussion) \citep{dall'anese, simonetto}.

Structured algorithms, on the other hand, are designed specifically for dynamic problems by incorporating a (possibly simplified) model of their time-variability.
The benefit of designing structured algorithms is that they generally achieve significantly lower, and in some scenarios zero, asymptotic error in tracking the solution trajectory \citep{simonetto}.
Different approaches to structured algorithm design have been proposed, with the use of control theoretical tools recently achieving great success for unconstrained \citep{Bastianello2024,bianchin}, constrained \citep{casti2023}, stochastic \citep{Casti2024,simonetto_nonlinear_2024}, and distributed \citep{vanweerelt2025} problems.

These algorithms typically leverage the internal model principle, meaning they 
rely on some prior knowledge of the dynamic modes that generate the time-varying signals involved in the online problem.
However, having such prior knowledge is oftentimes impractical \citep{bianchin}, or may be inaccurate in practice \citep{Bastianello2024}.

Therefore, this paper seeks to remove the need for prior knowledge of the internal model by using the tools of system identification. 
System identification provides a control-theoretic framework for extracting dynamic information directly from observed data 
\citep{Ljung1998}, and has been successfully applied in various fields including
control systems, signal processes, machine learning. Common methods include the least mean squares algorithm, Kalman filtering, and the recursive least squares (RLS) algorithm \citep{Guo1994}, the latter of which we leverage in this paper.
In an online optimization context then, the knowledge gained through the application of these techniques can be leveraged to construct an internal model and, consequently, develop a novel structured algorithm.

In this paper we address this objective, offering the following contributions:
\begin{itemize}
    \item We design a novel online algorithm which integrates the internal model-based design proposed in \citep{Bastianello2024} with a system identification procedure (in particular, recursive least squares) that serves to construct the internal model from observations of~\eqref{eq:objective}.

    \item We analyze the convergence of the resulting algorithm applied to online quadratic problems.

    \item We test the proposed algorithm for both quadratic and non-quadratic problems, showcasing its improved performance over the state of the art and, importantly, adaptability to changes in the internal model.
\end{itemize}

\section{Problem Formulation}\label{sec:problem}
In order to design the algorithm proposed in section~\ref{sec:algorithm}, we now focus our attention to a particular case of~\eqref{eq:objective}, characterized by the following quadratic function: 
\begin{align}\label{eq:quadfunc}
    f_k(\x):=\frac{1}{2}\x^\intercal\A\x +\x^\intercal \mathbold b_k,
\end{align}
where $\x\in \mathbb R^n$, the symmetric matrix $\A\in\mathbb R^{n\times n}$ is positive definite, and sequence $\{\mathbold b_k\in \mathbb R^n\}_{k\in \mathbb N}$ represents the time-varying signal driving the online problem.

Let $\lmin$ and $\lmax$ denote the minimal and maximal eigenvalues of matrix $\A$, respectively. Since $\A$ is positive definite, it follows that $\lmin\mathbold{I}\preceq\A =\A^\intercal\preceq \lmax\mathbold I$, and $\lmin>0$. Consequently, for any  $k\in\mathbb{N}$, each cost function $f_k$ is  $\lmin$-strongly convex and $\lmax$-smooth.
Moreover, it guarantees that, for any given sequence $\{\mathbold b_k\}$, each optimization problem in \eqref{eq:objective} admits a unique minimizer,  thereby defining a unique optimal trajectory $\{ \x_k^*= \A^{-1}\mathbold b_k \}_{k \in \N}$.
We remark that, despite the problem having a closed-form solution, we aim to design an algorithm that only uses gradient information. The goal indeed is for the algorithm to be applicable to a broader range of online problems, not only quadratic ones.

We introduce now the following model for $\boldsymbol{b}_k$ to enable our control-based algorithm design.

\begin{assumption}[Model of $\mathbold b_k$]\label{as:model-b}
We assume that the time series $\{\mathbold b_k\}_{k\in \mathbb N}$ admits a rational $\mathcal{Z}$-transform of the form 
\begin{align} \label{eq:Z-transform}
    \mathcal{Z}[\mathbold b_k] = \B(z) = \frac{\B_N(z)}{B_D(z)}, 
\end{align}
where real polynomials $B_D(z) = z^m+\sum_{i=0}^{m-1}d_iz^i$, and $\B_N(z) = \sum_{i=0}^p \boldsymbol{u}_i z^i$ with $p \leq m$. Furthermore, we assume that all roots of the denominator polynomial $B_D(z)$ have non-positive real parts.
\end{assumption}

Assumption~\ref{as:model-b} is standard in control-based design for online optimization problems, as it rules out unstable modes in $\mathbold b_k$, which would otherwise cause unbounded growth in the sequence of minimizers, \textit{i.e.} $\norm{\x_k^* - \x_{k-1}^*} \to \infty$ as $k \to \infty$.\
The difference with previous works such as \citep{Bastianello2024,casti2023}, is that \textit{we do not assume $B_D(z)$ (and $\B_N(z)$) to be known}, only that $\boldsymbol{b}_k$ admits the $\mathcal{Z}$-transform.

Finally, as discussed above we assume the algorithm has access to an \textit{oracle of the gradient}, as well as the bounds on $\A$'s eigenvalues, $\lmin$ and $\lmax$. This assumption is typical in structured online optimization algorithms
(see, \textit{e.g.}, \citep{Bastianello2024}).

\smallskip

In the setting discussed above, the goal then is to design a control-based algorithm which integrates online optimization
with a system identification routine, removing the need of prior knowledge on the evolution of $\mathbold b_k$. Applying identification indeed allows to asymptotically reconstruct the internal model (specifically, $B_D(z)$), and as a consequence achieving perfect tracking of the optimal trajectory.

\section{Proposed Algorithm}\label{sec:algorithm}
In this section we design the proposed algorithm, which we call \textit{SIMBO}, for \textit{Self-Identifying Internal Model-Based Online Optimizer}.
As discussed above, the foundation of our algorithm design is the algorithm proposed in \citep{Bastianello2024}, which is characterized by the updates:
\begin{subequations}\label{eq:CBALG}
\begin{align}
    \w_{k+1} &= (\mathbold F \otimes \mathbold I)\w_k + (\G \otimes \mathbold I) \nabla f_k(\x_k)\\
    \x_{k+1} &=(\mathbold K \otimes \mathbold I) \w_{k+1}
\end{align}
\end{subequations}
where $\w$ is the internal state of the algorithm, and
\begin{equation}\label{eq:state-space-realization}
\begin{split}
	\boldsymbol{F} &= \begin{bmatrix}
		0 & 1 & 0 & \cdots \\
		& & \ddots & \\
		0 & \cdots & 0 & 1 \\
		-d_0 & \cdots & \cdots & -d_{m-1}
	\end{bmatrix},\quad
	\G= \begin{bmatrix}
		0\\
		\vdots\\
		0\\
		1
	\end{bmatrix}\\
    \boldsymbol{K} &=\begin{bmatrix}c_0 & c_1 & \cdots & \cdots & c_{m-1}\end{bmatrix},
\end{split}
\end{equation}
with $\boldsymbol{K}$ that could be computed according to \citep{Bastianello2024} if $B_D(z)$ were known; but in the setting of this paper it is not.

Therefore, we propose to identify it, and using the output of the identification routine to characterize~\eqref{eq:CBALG}.
However, at initialization ($k = 0$) we do not have any historical data on the problem~\eqref{eq:objective} that could be used to identify the internal model beforehand. The idea then is to design a two-phase algorithm.
In the first phase, we apply online gradient descent (OGD) \citep{dall'anese} to the problem:
\begin{equation}\label{eq:ogd}
    \x_{k+1}=\x_k-h\nabla f_k(\x_k),
\end{equation}
where $h < 2 / \lmax$ is the step-size. Since OGD is unstructured, we can directly apply it to the problem without any model information. The output of OGD then can be used as data to feed the system identification routine, while also providing a (sub-optimal) sequence of decisions $\x_k$.
Once enough data on the problem has been collected to construct a first approximation of the internal model, we trigger the second phase. During the second phase we use the approximate internal model to construct $\boldsymbol{F}$ in~\eqref{eq:CBALG}, and to compute the controller $\boldsymbol{K}$. We can then switch to using the control-based algorithm~\eqref{eq:CBALG} to compute the decisions $\x_k$.
But we do not stop using system identification, since during the second phase we can refine the identified internal model or, more importantly, adapt to changes in the problem.

\smallskip

In section~\ref{subsec:sys-id} we discuss the system identification routine that we use during the two phases, and in section~\ref{subsec:algorithm} we lay out the overall proposed algorithm.

\subsection{System identification routine}\label{subsec:sys-id}
The system identification tool that we select to identify the internal model $B_D(z)$ is recursive least squares (RLS) \citep{Guo1994}.
Letting $\mathbold{d} = [d_0, \ldots, d_{m-1}]^\intercal \in \R^m$ be the vector of coefficients of the polynomial $B_D(z)$, then RLS identifies it using the recursion:
\begin{equation}\label{eq:RLS_basic}
    \hat{\mathbold {d}}_{k+1} = \hat{\mathbold{d}}_{k} + \mathbold L_k(\mathbold y_k-\mathbold{\phi}_k^\intercal\hat{\mathbold{d}}_k),
\end{equation} 
where $\mathbold y_k$ denotes the observation data, $\mathbold \phi_k$ denotes the regressor, and $\mathbold L_k$ denotes the gain \citep{Guo1994}:
\begin{equation}\label{eq:gain-rls}
    \mathbold L_k = \mathbold P_k\mathbold{\phi}_k (\alpha\mathbold I +\mathbold{\phi}_k^\top \mathbold P_k \mathbold{\phi}_k)^{-1},
\end{equation}
with
\begin{equation}
    \mathbold P_{k+1}=\frac{1}{\alpha}\left(\mathbold P_k -\mathbold P_k\mathbold{\phi}_k (\alpha\mathbold I +\mathbold{\phi}_k^\top \mathbold P_k \mathbold{\phi}_k)^{-1}\mathbold{\phi}_k^\top\mathbold P_k\right),
\end{equation} 
where $\mathbold P_0 > 0$, and $\alpha \in (0,1)$ is a forgetting factor which weights recent entries more heavily \citep{Guo1994}.
The following paragraphs delineate the specific expressions that $\mathbold y_k$ and $\mathbold \phi_k$ take during the two phases of the algorithm when either~\eqref{eq:ogd} or~\eqref{eq:CBALG} are applied.

\subsubsection{Phase 1: initializing the identification}
During the first phase starting at $k = 0$ we apply the online gradient descent characterized by~\eqref{eq:ogd}. Therefore, we need to extract information about $B_D(z)$ from OGD as follows.
Taking the $\mathcal{Z}$-transform of the output of OGD, $\mathcal{Z}\{\x_{k+1}\} = \mathcal{Z}\{\mathbold x_{k} - h \nabla f_k(\x_k) \}$ and using~\eqref{eq:quadfunc}, yields
\begin{equation}\label{eq:timeshiftz}
    z\mathbold X(z) = (\mathbold I - h \mathbold A)\mathbold X(z) - h \mathbold B(z).
\end{equation}
Rewriting the $\mathcal{Z}$-transforms as infinite sums \citep{Graf2004}, and recalling Assumption~\ref{as:model-b}, we get
\begin{align*}
    &\sum_{k=0}^\infty \left( \mathbold x_{k+m}  +\sum_{i=0}^{m-1}d_i\mathbold x_{k+i}\right)z^{-k} = \\ &\hspace{3cm} -h \sum_{k=0}^\infty(\mathbold I-h \mathbold A)^kz^{-k-1} \sum_{j=0}^p \mathbold u_jz^j
\end{align*}
which then gives
\begin{equation}\label{eq:recurrugly1}
    \mathbold x_{k+m}  +\sum_{i=0}^{m-1}d_i\mathbold x_{k+i} = -h \sum_{j=0}^p(\mathbold I - h \mathbold A)^{k-1+j}\mathbold{u}_jz^j.
\end{equation}
Since we have selected a step-size $h < 2 / \lmax$ then we know that:
\begin{equation}
    \lim_{k\to\infty}\left(\mathbold I - h \mathbold A\right)^{k-1+j} = \boldsymbol 0,
\end{equation}
and thus the recurrence~\eqref{eq:recurrugly1} becomes: \begin{align} \label{eq:recurrugly2}
        \mathbold x_{k+m} + \sum_{i=0}^{m-1}d_i\mathbold x_{k+i} = \boldsymbol 0, \hspace{2ex} \forall k \geq m +1.
\end{align}
With this recurrence in place, we can then apply RLS with $\mathbold y_k = \mathbold{\phi}_k^\top\mathbold{d}$ and
\begin{equation}\label{eq:regressors}
    \mathbold{\phi}_k = [-\mathbold x_{k-1}, -\mathbold x_{k-2}, \cdots , -\mathbold x_{k-m}]^\top,
\end{equation}
and $\mathbold L_k$ as introduced in~\eqref{eq:gain-rls}.

\subsubsection{Phase 2: continual identification}
Once phase 1 has successfully constructed an approximation of the internal model, characterized by $\hat{\mathbold {d}}_k$, we can use it to define $\mathbold{F}$ according to~\eqref{eq:state-space-realization}, and to compute the controller $\mathbold{K}$.
We can then deploy the control-based algorithm~\eqref{eq:CBALG} to replace OGD.

However, the first phase only approximates the internal model, and in addition the model might change over time. Therefore, during the second phase we continue running RLS to improve the identified model, and to do so, we need to extract information from the output of~\eqref{eq:CBALG} as follows.
Letting the inexact model computed at the end of phase 1 (iteration $k_1$) be
\begin{equation*}
    \hat{B}_D(z) = z^m+\sum_{i=0}^{m-1} \hat{d}_{k_1, i} z^i,
\end{equation*}
and using the controller $C(z) = {C_N(z)} / {\hat{B}_D(z)}$, algorithm~\eqref{eq:CBALG} is represented by the $\mathcal{Z}$-transform:
\begin{equation*}
    \mathbold X(z) = \left(\hat{B}_D(z) \mathbold I-C_N(z)\mathbold A\right)^{-1}\frac{\mathbold B_N(z)C_N(z)}{B_D(z)}
\end{equation*}
and rearranging:
\begin{equation}\label{eq:CB-rec-der}
    B_D(z)\mathbold X(z)  = \left(\hat{B}_D(z) \mathbold I-C_N(z)\mathbold A\right)^{-1}\mathbold B_N(z)C_N(z).
\end{equation}
Again using the properties of the $\mathcal{Z}$-transform we rewrite~\eqref{eq:CB-rec-der} as
\begin{equation}
\begin{split}
\sum_{k=0}^\infty \Big(\mathbold x_{k+m}
  +\sum_{i=0}^{m-1} d_i \mathbold x_{k+i}\Big) z^{-k} = \\ \qquad \Big(\hat{B}_D(z)\mathbold I - C_N(z)\mathbold A\Big)^{-1} \mathbold B_N(z)\,C_N(z).
\end{split}
\end{equation}
Finally, since the right-hand side is made up of two anticausal signals, and $C_N(z)$ is chosen in such a way that $\left(\hat{B}_D(z) \mathbold I-C_N(z)\mathbold A\right)$ is stable, the recurrence $\mathbold x_{k+m}  +\sum_{i=0}^{m-1}d_i\mathbold x_{k+i} = \boldsymbol 0$, $k\geq m+1$, holds (same as~\eqref{eq:recurrugly2}).
We can then apply the RLS~\eqref{eq:RLS_basic} with $\mathbold y_k = \mathbold{\phi}_k^\top\mathbold{d}_k$ and $\mathbold{\phi}_k = [-\mathbold x_{k-1}, -\mathbold x_{k-2}, \cdots , -\mathbold x_{k-m}]^\top $.

\subsection{Algorithm overview}\label{subsec:algorithm}
\begin{figure}[!ht]
    \centering
    \includegraphics[width=\linewidth]{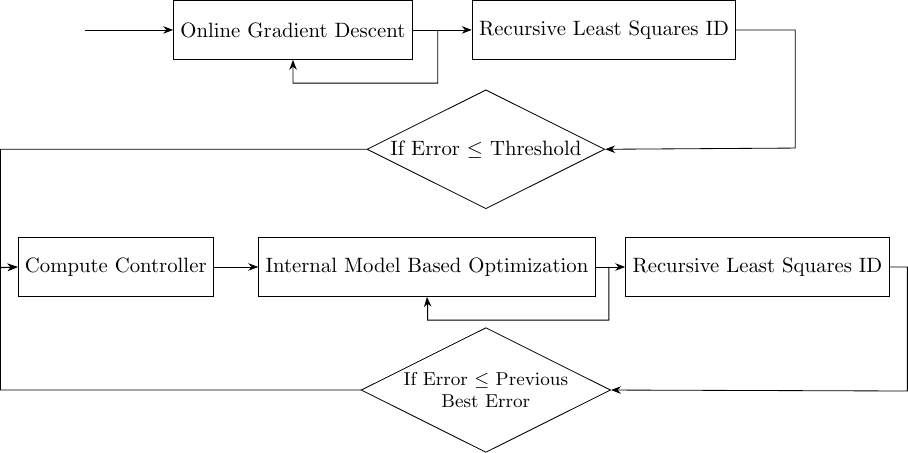}
\caption{Flowchart of the SIMBO algorithm}
\label{fig:alg_form}
\end{figure}

We are now ready to formalize the overall algorithm that we propose, SIMBO.
Figure~\ref{fig:alg_form} represents the flowchart of SIMBO, highlighting the two phases, 1. identification initialization with online gradient descent, and 2. continual identification with the control-based algorithm~\eqref{eq:CBALG}.
We now need to define the triggering condition that switches from phase 1 to 2, and the condition that triggers a recomputation of the controller in phase 2.

Given that RLS is applied to the recurrence~\eqref{eq:recurrugly2}, we can evaluate the system identification performance by the error $e_k := \norm{\mathbold y_k-\mathbold{\phi}_k^\intercal\hat{\mathbold{d}}_k}_1$. If the estimated coefficients $\hat{\mathbold{d}}_k$ coincide with the coefficients of the actual model this error is zero.
Therefore, we can select a threshold $\theta > 0$ such that, if $e_k \leq \theta$ then we know that phase 1 has constructed a good approximation of the internal model, and we can reliably use this model to compute the controller characterizing~\eqref{eq:CBALG} and switch to phase 2.
The trigger in phase 2, instead, is to ensure that we continuously improve the approximation of the internal model, especially in the case of changes of the actual model. This trigger then determines a recomputation of the controller when $e_k \leq e_{\bar{k}}$ where $e_{\bar{k}}$ was the previous best error. In other words, this condition is satisfied when the approximation of the model has improved.

\subsubsection{Practical heuristics}
To complete the algorithm, we now discuss two heuristics that improve the performance in practice.
First of all, in principle the triggering condition during phase 2 might be verified at each iteration $k$, thus requiring a large number of controller recomputations, which require the solution of two LMIs \citep{Bastianello2024}. This, however, would increase significantly the computational complexity and the risk of incurring in infeasible controller design problems. The idea then is to allow a recomputation of the controller only if the identification error $e_k$ has not improved for a number of iterations; that is $e_k \leq e_{\bar{k}}$ and $k \geq \bar{k} + t$.
Additionally, if the controller design problem happens to be infeasible we easily fall back on the previous controller.

The second heuristic we integrate in the algorithm is how to deal with changes in the actual internal model. The idea is to trigger \textit{phase 1} again once the identification error has worsened significantly. In particular, we trigger phase 1 when
$$
    \norm{\mathbold y_k-\mathbold{\phi}_k^\intercal\mathbold{d}_{\bar{k}}}_1 > C \norm{\mathbold y_{k-1}-\mathbold{\phi}_{k-1}^\intercal\mathbold{d}_{\bar{k}}}_1,
$$
with $C \gg 1$ and where $\mathbold{d}_{\bar{k}}$ is the identified model with the best error $\bar{k}$ up to time $k-1$ (in section~\ref{sec:Numerics} we use $C = 100$).

\section{Convergence Analysis}\label{sec:convergence}
In this section we discuss the convergence of the proposed algorithm SIMBO. To this end, we assume that the actual internal model does not change (\textit{i.e.} Assumption~\ref{as:model-b} holds). This means that phase 1 is executed once, and then we switch to phase 2 (no heuristics are applied).

\begin{proposition}
Let Assumption~\ref{as:model-b} hold, and assume that $\{b_{k}\}_{k\in \mathbb{N}}$ is persistently exciting of order $m$
\footnote{The real-valued sequence $\{b_{k}\}_{k\in \mathbb{N}}$ is said to
be persistently exciting (PE) of order $m$, if the Hankel matrix
$H_m(b_{[0, h-1]})$ has full row rank for some integer $h\geq m$. Here $H_m(b_{[0, h-1]})=[\psi_{0}, \psi_{1}\cdots, \psi_{h-m}],$ where $\psi_{k}=[b_{k}, \cdots, b_{k+m-1}]^{\top} (0\leq k\leq h-m).$ \citep{Willems_2005}.}.
Then the output of SIMBO, $\{ \x_k \}_{k \in \N}$ verifies
$$
    \lim_{k \to \infty} \norm{\x_k - \x_k^*} = 0.
$$
\end{proposition}
\paragraph*{Proof}
We start by analyzing the convergence during phase 1. The data used in the RLS comes from the online gradient descent~\eqref{eq:ogd}, and we need to show that $\hat{\mathbold{d}}_k$ of~\eqref{eq:RLS_basic} is indeed converging towards the internal model coefficients $\boldsymbol{d}$.
First of all, we remark that OGD is converging to a neighborhood of the optimal trajectory, as proved in the following. By \cite[Theorem~1]{simonetto} we have
\begin{equation}\label{eq:bound-ogd}
    \norm{\mathbold x_{k+1} -\mathbold x_{k+1}^* } \leq\varrho\norm{\mathbold x_k-\mathbold x_k^*}+\norm{\mathbold x_k^* - \mathbold x_{k+1}^*}
\end{equation}
where $\varrho := \max \{ |1 - h \lmin|, |1 - h \lmax| \} \in (0,1)$ for $h < 2/\lmax$, and where
$$
    \norm{\mathbold x_k^* - \mathbold x_{k+1}^*} \leq\norm{\mathbold A^{-1}} \norm{\mathbold b_k - \mathbold b_{k+1}}\leq\frac{1}{\underline{\lambda}}\norm{\mathbold b_k - \mathbold b_{k+1}}
$$
since $\x_k^* = - \boldsymbol{A}^{-1} \boldsymbol{b}_k$.
Let $k_1 \in \N$ be the time when the algorithm switches to phase 2; then by~\eqref{eq:bound-ogd} we have
\begin{align*}
    \norm{\x_{k_1} - \x_{k_1}^*} &\leq \varrho^{k_1} \norm{x_0 - \x_0^*} \\ &+ \frac{1 - \varrho^{k_1+1}}{1 - \varrho} \max_{k \in [0, k_1]} \frac{1}{\underline{\lambda}}\norm{\mathbold b_k - \mathbold b_{k+1}}
\end{align*}
and OGD indeed converges to a neighborhood of the optimal trajectory during phase 1.

Now, by~\eqref{eq:ogd} and the fact that the cost is quadratic, the output of~\eqref{eq:ogd} can be written as 
\begin{equation}\label{eq:recursion-x}
        \mathbold x_{k+1} = (\mathbold I -h \mathbold A)^{k+1} \mathbold x_0-h\sum_{j=0}^k(\mathbold I - h \mathbold A)^{k-j}\mathbold b_j,
\end{equation}
which generates the regressors $\boldsymbol{\phi}_k$ according to~\eqref{eq:regressors}.
The RLS then is converging provided that, defining
$
    \mathbold\Phi_{k}=\begin{bmatrix}
    \mathbold \phi_{k+m} & \cdots & \mathbold \phi_{k+h-m} \end{bmatrix},
$
the matrix $\mathbold\Phi_k \mathbold \Phi_k ^\intercal$ is invertible for each $k\geq 0$.
But this is a consequence of the fact that there is no noise in~\eqref{eq:recursion-x}, that $\{b_k\}_{k\in \mathbb{N}}$ is persistently exciting of order $m$, and that $\boldsymbol{I} - h \boldsymbol{A}$ is full rank for any $h < 2 / \lmax$ \citep{haykin}.
In this scenario, we can characterize that the regressor sequence $\{x_{k}\}_{k \in \N}$ is persistently exciting of order $m$, which yields the following bound on the identification error \citep{Guo1994, haykin}: $\norm{\hat{\mathbold{d}}_{k} - \mathbold{d}} \leq M\zeta^{k} \norm{\hat{\mathbold{d}}_{0} - \mathbold{d}}$ for some $\zeta \in (0, 1)$ and $M>0$.
This implies that indeed, at the end of phase 1, RLS has approximately identified the internal model. Thus, there exists $k_1 \in \N$ large enough so that phase 2 is triggered.

\smallskip

We are now ready to analyze the convergence of phase two in $(k_1, \infty)$.
During this phase, we concurrently run RLS, which outputs $\hat{\mathbold{d}}_k$, and the control-based algorithm~\eqref{eq:CBALG}, whose matrices $\mathbold{F}_k$ and $\mathbold{K}_k$ are computed using $\hat{\mathbold{d}}_k$.
By \citep[Proposition 4]{Bastianello2024} we know that~\eqref{eq:CBALG} converges to a bounded neighborhood of the optimal trajectory when it is using an inexact internal model $\hat{\mathbold{d}}_k$, and that the neighborhood shrinks as $\hat{\mathbold{d}}_k \to \mathbold{d}$.
To ensure this happens, we need then to guarantee that the RLS is converging to $\mathbold{d}$. This, similarly to the proof for phase 1, is a consequence of $\mathbold{b}_k$ being persistently exciting. This is because when $\hat{\mathbold{d}}_k$ closely matches $\mathbold{d}$, the regressor is determined by the the output of the internal model-based algorithm of \citep{Bastianello2024}. Since we know that the asymptotic tracking error of this algorithm is zero from \citep[Proposition 1]{Bastianello2024}, the output of the algorithm is approximately $\mathbold x_{k+1}\approx\mathbold A^{-1}\mathbold b_k$. Hence the matrix $\mathbold{\Phi}$ is defined by $\mathbold b_k$ like in phase 1 as the matrix $\boldsymbol A$ is not time varying. The convergence then depends on the persistent excitation of $\mathbold{b}_k$. $\hfill\square$

\section{Numerical Results}\label{sec:Numerics}
In this section we analyze the performance of SIMBO, and compare it to OGD~\eqref{eq:ogd} and the control-based algorithm~\eqref{eq:CBALG} \citep{Bastianello2024}. The latter algorithm is designed using the exact internal model; as such, it serves as a baseline for SIMBO, but as discussed before the exact model would not be available in practice.
All simulations were implemented using the \texttt{tvopt} Python package \citep{tvopt}.

\subsection{Quadratic problems}\label{subsec:quadratic-numerical}
We start comparing the three algorithms for quadratic problems characterized by~\eqref{eq:quadfunc}, with $n = 15$, and $\boldsymbol{A}$ is randomly generated so that $\lmin = 1$, $\lmax = 5$.
For the linear term $\boldsymbol{b}_k$ we use the following four internal models:
\begin{enumerate}
    \item $\mathbold b_k = \sin(\omega_0k T_s)\boldsymbol 1$ \quad  (Sine function)
    \item $\mathbold b_k = k T_s\bar {\mathbold b}$ \quad (Ramp function)
    \item $\mathbold b_k = \sin(\omega_0k T_s)\boldsymbol 1 +  k T_s\bar {\mathbold b}$ \quad (Sine plus ramp function)
    \item $\mathbold b_k = \sin^2(\omega_1k T_s)\boldsymbol 1$  \quad (Sine squared function)
\end{enumerate}
where $T_s = 0.1$, $\omega_0 = 1$, $\omega_1 = 10$, and $\bar{\mathbold b}\in\mathbb{R}^n$ is randomly generated.

\begin{figure}[!ht]
    \centering
    \includegraphics[width=\linewidth]{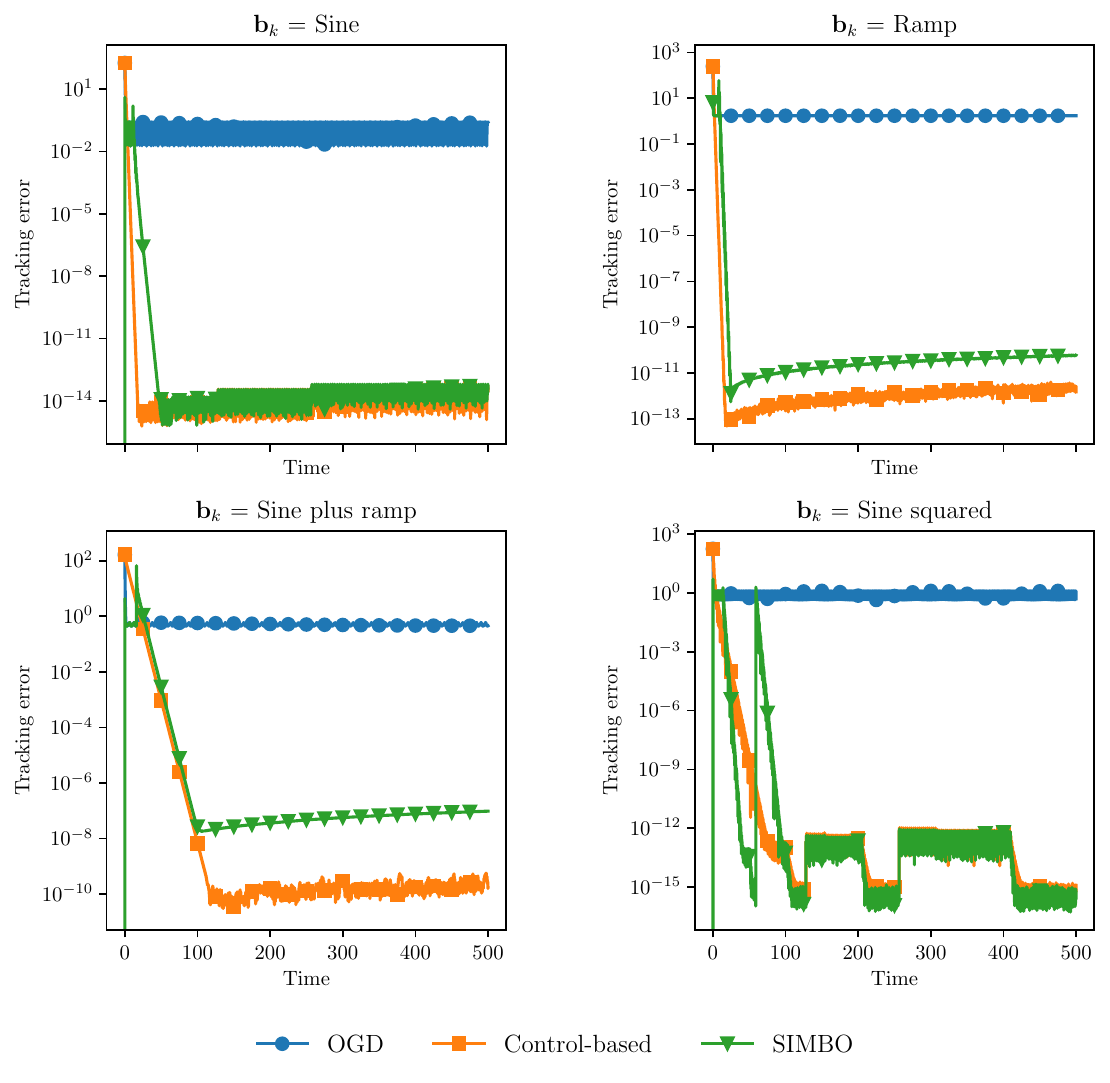}
    \caption{Tracking error comparison with different internal models in~\eqref{eq:quadfunc}}
    \label{fig:TE_1}
\end{figure}

In Figure~\ref{fig:TE_1} we report the evolution of the tracking error $\{ \norm{\x_k - \x_k^*} \}_{k \in \N}$ for the three algorithms. 
As expected, OGD can only converge to a bounded neighborhood of the optimal trajectory, while the control-based algorithm converges asymptotically to it (up to numerical precision). This is of course dependent on knowing the exact model, which in practice might not be available. Therefore, seeing that SIMBO performs very closely to the control-based algorithm guarantees that indeed incorporating system identification allows to reduce the prior knowledge requirements without sacrificing accuracy. In other words, the system identification routine successfully identifies the internal model.
Inspecting Figure~\ref{fig:TE_1}, some additional observations are in order.
First of all, the cases where a ramp signal is involved, numerical precision is slightly worse since the ramp grows unbounded. Incorporating an integrator in the algorithm design might serve to reduce this numerical issue.
The second observation is that switching from phase 1 to 2, or triggering the recomputation of the controller (bottom right plot at around $k = 50$), might give rise to a transient. Further changes to the design might serve to reduce the size of these transients.

We conclude this section by providing the asymptotic value of the tracking errors in Table~\ref{tab:ATE_table}. The asymptotic error is estimated in practice by taking the maximum error over the final $4/5$ of the simulation.
\begin{table}[!ht]
\caption{Asymptotic tracking errors}
\label{tab:ATE_table}
    \centering
    \begin{tabular}{c c c c c}
         Algorithm & Ramp & Sine & Sine$^2$ & Sine + ramp \\
        \hline
        OGD  & $1.73$e$+00$ & $2.70$e$-01$ & $1.29$e$+00$ & $5.97$e$-01$ \\
        Control-Based & $4.02$e$-12$& $5.51$e$-14$& $8.03$e$-13$& $5.65$e$-10$  \\
        SIMBO & $6.09$e$-11$& $6.40$e$-14$& $7.28$e$-13$& $1.90$e$-07$  \\
       \hline
    \end{tabular}
\end{table}

\subsection{Adapting to a changing internal model}
\begin{figure}[!ht]
    \centering
    \includegraphics[width=\linewidth]{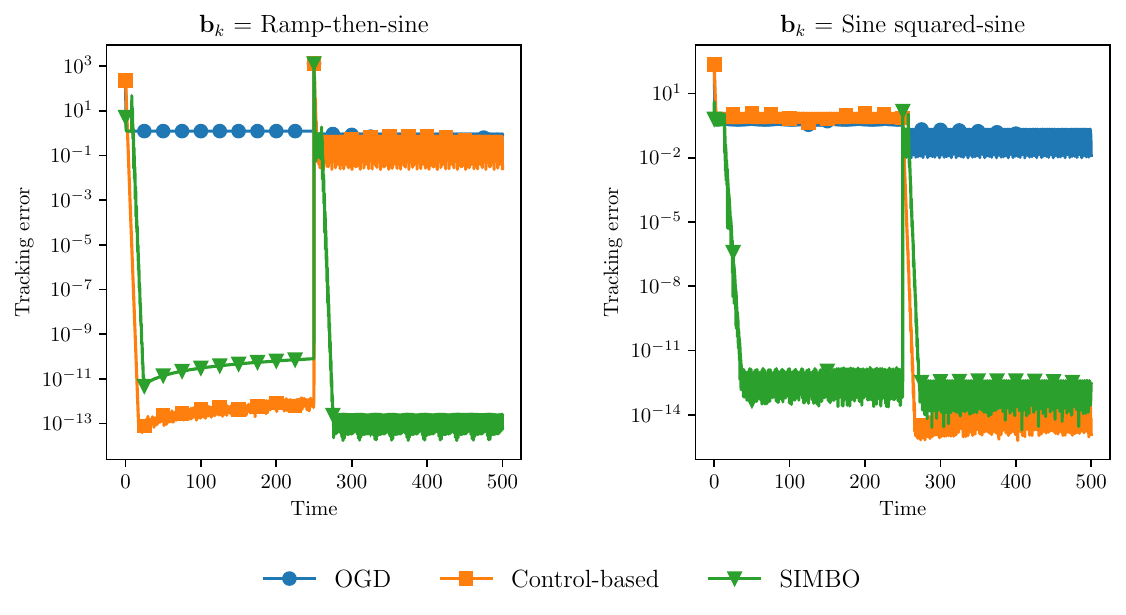}
    \caption{Tracking error with changing internal models}
    \label{fig:TE_2}
\end{figure}

The use of SIMBO, as opposed to the control-based algorithms~\eqref{eq:CBALG}, is especially necessary when the internal model might change over time.
In this section we test the algorithms for two quadratic problems whose linear term $\boldsymbol{b}_k$ changes internal model half-way through, (1) ramp then sine, (2) sine then squared sine (using the models of section~\ref{subsec:quadratic-numerical}).
The results are depicted in Figure~\ref{fig:TE_2}.
We can see that, as expected, SIMBO is able to adapt to changes in the internal model. Indeed, after the change, SIMBO re-identifies the model (by triggering phase 1 again, see section~\ref{subsec:algorithm}), recomputes the controller, and then successfully tracks the optimal trajectory.
On the other hand, the control-based algorithm is fixed, based on the internal model used to compute the controller at the beginning. Thus, when the internal model coincides with the model of $\boldsymbol{b}_k$ (first half in the left plot, second half in the left plot) the algorithm tracks the optimal trajectory. However, as soon as the model changes half-way through the simulation, its performance decays, in some cases worse than the unstructured OGD.

\subsection{Time-varying Hessians}
In this paper we have designed SIMBO for problems with quadratic cost~\eqref{eq:quadfunc}, where only the linear term varies. However, as proved in \citep{Bastianello2024}, this control-theoretical approach can equally be applied to more general problems, and the same applies to SIMBO\footnote{Indeed, the algorithm only requires oracle evaluations of the gradient to be applied in practice.}.
In this section then we test the algorithms on the following costs where also the Hessian changes over time:
$
    f_k(\mathbold x) = \frac{1}{2}\mathbold x^\intercal \mathbold A_k \mathbold x + \mathbold x ^\intercal \mathbold b_k,
$
where $\mathbold A_k =  \mathbold A + \tilde{\mathbold A}_k$, $\mathbold A = \mathbold V\mathbold \Lambda\mathbold V^\intercal$ and $\tilde{\mathbold A}_k =\mathbold V\text{diag}\{\sin(\omega_0kt_s)\mathbold v\}\mathbold V^\intercal$. Similarly to section~\ref{subsec:quadratic-numerical}, we guarantee that $\mathbold A_k$ has eigenvalues in $[1, 5]$ for all $k \in \N$.
For the sake of simplicity, we also assume that $\mathbold b_k = \bar{\mathbold b}\in\mathbb{R}^n$.

\begin{figure}[!ht]
    \centering
    \includegraphics[width=\linewidth]{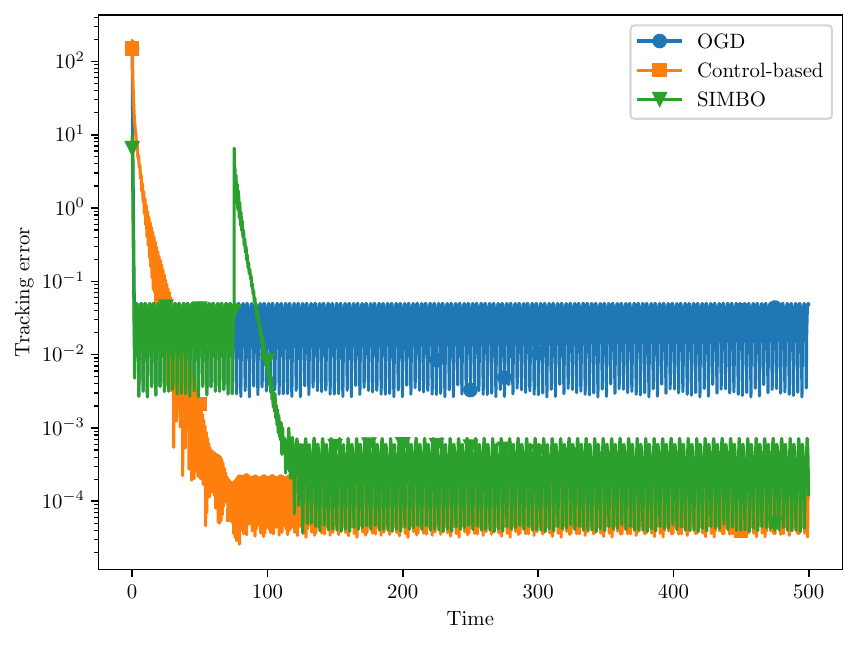}
    \caption{Tracking error with time-varying Hessian}
    \label{fig:TV_hessian}
\end{figure}

Figure~\ref{fig:TV_hessian} reports the evolution of the tracking errors in this setting.
As expected, OGD reaches a large neighborhood of the optimal trajectory, while the control-based algorithm~\eqref{eq:CBALG} and SIMBO reach a tighter neighborhood (although not zero, as the assumptions of section~\ref{sec:problem} are not verified).
Importantly, the control-based algorithm~\eqref{eq:CBALG} is defined on an internal model that is hand-tuned to improve performance (which requires prior information), while SIMBO reaches similar performance automatically adapting the internal model.

\section{Conclusions and Future Work}\label{sec:conclusions}
In conclusion, this paper proposes a novel class of structured online algorithms that merge the control theoretical design of \citep{Bastianello2024} with a system identification routine. The use of system identification allows to construct an internal model of the problem's time-variability, without having to resort to prior knowledge/data which in practice are rarely available. Additionally, using identification allows to adapt in real time to changes in the behavior (\textit{i.e.} in the internal model) of the online problem. The performance of the proposed algorithm is evaluated theoretically and validated numerically.
Future research will focus for example on extending this identification approach to more general problems (using the non-linear internal model principle).

\bibliography{ifacconf}
\end{document}